\documentclass[twoside,11pt]{article}

\usepackage{jmlr2e}

\usepackage{algorithm}% http://ctan.org/pkg/algorithms
\usepackage{algpseudocode}% http://ctan.org/pkg/algorithmicx
\usepackage{amsmath}
\usepackage{graphicx}
\usepackage{bbm}
\usepackage{tabulary}
\usepackage{booktabs}
\usepackage[symbol]{footmisc}

\jmlrheading{1}{2020}{1-48}{4/00}{10/00}{Larkin Liu}

\ShortHeadings{Algorithm for Two-Phase Facility Planning}{Liu 2020}
\firstpageno{1}

\begin{document}

\title{Algorithm for Two-Phase Facility Planning via Balanced Clustering and Integer Programming}

\author{\name Larkin Liu \email larkin.liu@mail.utoronto.ca \\
       \addr University of Toronto
}

\maketitle

\begin{abstract}
We present a solution for a two-phase facility planning scenario where in the first phase, there is some flexibility in determining where the locations of facilities (or sources) should fall. And in the second phase, new waypoints (or sinks) are added, but the location of the facilities are static. This solution applies the use of balanced clustering - using a modified K-Means approach, ensuring the cardinality of each group to be equal. Subsequently, it is followed by an integer programming solution, to solve the Hitchcock Transportation Problem. We show that the final solution can justifiably approximate the near optimal solution, and be a successful guide for facility planning in this specific scenario.
 \\
\end{abstract}

\begin{keywords}
  Balanced Clustering, Integer Programming, Facility Planning
\end{keywords}

\section{Introduction}

The optimization of facilities - also referred to as \textit{sources} in the network flow formulation - in modern supply chains play a critical role in the profitability and market competitiveness of any modern business. Such facilities can represent the function of a distribution center for goods or supplies; in our work we refer to these as depots. Such distribution centers can supply the stores or service centers in the surrounding location. In the network flow formulation these are referred to as \textit{sinks}, in our work we refer to them as waypoints.

Therefore, it is essential to ensure that in the path planning of any supply chain network is properly optimized to reduce transportation distance. We wish to find an appropriate solution for a unique problem that involves the positioning and assignment of depots and waypoints in a realistic situation. We are presented with a two-phase optimization process, in the first phase, the Positioning Phase denoted as Phase I, we must determine where facilities must fall on an unconstrained $\mathbb{R}^2$ plane, to obtain a fixed positioning of depots. Subsequently we are presented with the Assignment Phase denoted as Phase II, where an arbitrary set of waypoints are placed in $\mathbb{R}^2$ representing new locations of waypoints.

For the new and preexisting waypoints the distance from itself to each depot is known, and thus is represented by a set of coordinates are transformed into a symmetric distance matrix, or cost-matrix, $\mathbb{M}$. This cost-matrix is based the Euclidean distance from any waypoint to any depot. From this we must assign the optimal waypoint to depot assignment such that the total travel distance is minimized. 

The situation stipulates that certain locations, or waypoints, must be assigned to respective facilities, or depots. Phase I presents an optimization problem where the location of the depot lies unconstrained in $\mathbbm{R}^2$, thus we option to use a clustering algorithm to solve iteratively, specifically the K-Meams algorithm.

The K-Means algorithm can also be interpreted as the Expectation Maximization (EM) algorithm \citep{Dempster:1977}, without the use of a probability distribution. Specifically the assignment step involves assigning each waypoint to its nearest neighbour depot. The result provides a series of Voronoi cells \citep{Voronoi:1908} in $\mathbbm{R}^2$, where each depot forms the centroid of each Voronoi cell. Yet, we understand that the major limitation of K-Means for our use-case is that it does not ensure that the number of waypoints belonging to each Voronoi cell is equal. For this, we must apply a modified K-Means implementation authored by \citep{Schubert:2015}, ensuring that each cluster of waypoints contains the same number, and the centroid of each cluster is the location of the depot. We refer to this as K-Means Clustering with Equal Size Constraints, or balanced clustering.

In Phase II, the situation is fundamentally different. There exists no freedom in the movement of the depots in $\mathbbm{R}^2$, and thus the problem is purely an assignment problem, which waypoints to assign to with depots. All waypoints, including the preexisting ones from Phase I as well as newly added ones from Phase II undergo a reassignment.

The solution to this type of problem has been investigated for example in \citep{Hakimi:1964} offering a reformulation of the Weber Problem \citep{Weber:1909} for facilities, where each facility in $\mathbbm{R}^2$ can function as both a depot or waypoint. \citep{Hakimi:1964} proves that the optimal solution must intersect with the set of facilities. In other words, the optimal depots locations correspond to facility locations. However, in our application there is a distinction in the type of facility, that is depots and waypoints are not interchangeable. Thus we must model the problem as an Hitchcock Transportation Problem \citep{Hitchcock:1941}, or HTP. The HTP is also another extension of the Weber Problem where, the placement of each depot and waypoint are fixed, non-interchangeable, and a solution must be provided to assign each waypoint to a depot, such that some distance function is reduced.

In the business setting, the introduction of new waypoints can represent a decision to add more location to a store branch, as determined external sources, such as executive management, at a later time. The idea is that a full-disclosure of where the waypoints location are, and how many are to be planned, is not given ahead of time.  For example, an initial deployment of store branches followed by a second wave of store openings. Where in the initial phase of optimal locations must be planned for the supply depots in terms of minimizing travel distance. Therefore, there is uncertainty regarding where or how many of the waypoints can be placed. We only have control over the placement of the depots in the first phase, and assignment of waypoints to depots in the second phase.

\subsection{Experimental Outline} \label{sec:data_desc}

The data content we select to model consists of approximately 25,000 world wide locations of Starbucks chain store locations. The store location's longitude and latitude are provided. We assume that all distances from any location on the 2D world map to be the Euclidean distance, with no geographical, topographical, land, or ocean boundary which can affect the Euclidean distance. We also assume, that from the initial placement of waypoints in Phase I, there exists no correlation or deterministic relationship to the subsequent waypoint placements in Phase II.  

We begin by selection at random percentage of waypoints, represented by $\gamma = 0.05$, for Phase I. It is evident that the algorithm presented in Phase II will solve for the majority of the total waypoint to depot assignment. Subsequently, we introduce the rest of the waypoints, $1-\gamma$, for the assignment phase. In each of the respective phases, we present a specific solution to the different scenarios. In Phase I we apply the balanced clustering algorithm \citep{Schubert:2015}, and in Phase II we present a solution to the HTP. Both to be illustrated with further detail in the later sections.

\begin{table}[htbp]\caption{Summary of Notation}
\begin{center}% used the environment to augment the vertical space
% between the caption and the table
\begin{tabular}{r c p{10cm} }
\toprule
$\mathbb{D}$ & $\triangleq $ & Set of all potential depots.\\
$\widetilde{\mathbbm{D}}$ & $\triangleq$ & Set of all depots after Phase I.\\
$d$ & $\triangleq$ & Location of an individual depot.\\
$d^*$ & $\triangleq$ & Optimal location of individual depot for Phase I.\\
$d'$ & $\triangleq$ & Best alternate depot assignment for a waypoint $w$ corresponding to $\Delta$.\\
$\mathbb{M}$ & $\triangleq$ & Distance matrix representing distance from all possible waypoints to depots.\\
$N$ & $\triangleq$ & Total number of waypoints in either Phase I or Phase II.\\
$N_I$ & $\triangleq$ & Total number of waypoints in Phase I.\\
$N_{II}$ & $\triangleq$ & Total number of waypoints in Phase II.\\
$n_k$ & $\triangleq$ & Number of waypoints assigned to each depot.\\
$K$ & $\triangleq$ & Number of depots.\\

$\mathbb{W}$ & $\triangleq$ & Set of all possible waypoints.\\
$\widetilde{\mathbbm{W}}$ & $\triangleq$ & Set of new waypoints introduced in Phase II.\\
$\widehat{\mathbbm{W}}$ & $\triangleq$ & Set of all waypoints in Phase I and II.\\
$\overline{\mathbbm{W}}$ & $\triangleq$ & Sorted $\widehat{\mathbbm{W}}$ based on Euclidean distance difference between $w$ to $d$ versus to $d'$\\
$\overline{\mathbbm{W}}_d$ & $\triangleq$ & Subset of $\overline{\mathbbm{W}}$ corresponding to all $w$ assigned to $d$.\\
$\overline{\mathbbm{W}}_{d'}$ & $\triangleq$ & Subset of $\overline{\mathbbm{W}}$ corresponding to all $w$ assigned to $d'$.\\

$w$ & $\triangleq$ & An individual waypoint.\\
$w$ & $\triangleq$ & An individual waypoint.\\

$\Delta$ & $\triangleq$ & The \textit{assignment plan}, defined as the set of all tuples (w, d).\\
$\Delta^*$ & $\triangleq$ & The optimal \textit{assignment plan}, as defined in Eq. (\ref{eq:opt_def}).\\
$\Delta_d$ & $\triangleq$ & Subset of $\Delta$ representing the set of waypoints only assigned to depot $d$.\\
$\Delta'$ & $\triangleq$ & Assignment plan corresponding to the assignment of all waypoints in $\Delta$ to its second closest alternate depot.\\

$\Gamma$ & $\triangleq$ & Set of all possible distances from all waypoints to all depots.\\
$\psi(w, d)$ & $\triangleq$ & Distance function of waypoint to depot.\\
$\Psi(\Delta)$ & $\triangleq$ & Total cost of the assignment plan in terms of Euclidean distance.\\

$\mathbbm{1}_{\Delta}$& $\triangleq$ & \(\left\{\begin{array}{rl}
1,  & \text{if $w$ is assigned to $d$.} \\
0,  & \text{otherwise} \end{array} \right.\)\\

\bottomrule
\end{tabular}
\end{center}
\label{tab:notation}
\end{table}

\section{Problem Statement} \label{sec:prob_statement}

The number of depots, $K$, are fixed and equal for both Phase I and Phase II. The waypoints are given in two phases, we let $N_I$ represent the number of waypoints in the Phase I and $N_{II}$ for Phase II respectively. In our notation, $N$ flexibly represents the number of waypoints in either Phase I or Phase II. Most importantly, the number of waypoints assigned to each depot must be equal, outlined in Eq. \ref{eq:nk}. In our specific problem, we pose a situation where only a subset of the entire set waypoints are initially know. 

\begin{equation} \label{eq:nk}
n_k = N/K 
\end{equation}

In the Phase I, or the \textit{Positioning Phase} some waypoints are provided, however, the depot locations lie unconstrained in $\mathbb{R}^2$, and can be modified accordingly. We can infer that the distance function is not discrete, rather can be modeled as a continuous function. In this experiment we use the simple Euclidean distance in $\mathbb{R}^2$, between two 2-dimensional points. The earliest notion of such an unconstrained facility planning problem can be illustrated in the Weber Problem \citep{Weber:1909}. Its objective is to find a single depot by minimizing the Euclidean distance between all waypoints and its nearest depot. This problem can be solved using various geometric methods, such as some illustrated in \citep{Fekete:2003}. The big drawback using such geometric methods are that they are not easily extendable beyond $\mathbb{R}^2$, and is not easily programmable into computers. For this reason, we propose to use computational methods, such as balanced clustering, to solve Phase I.

After we present Phase II, or the \textit{Assignment Phase}, where additional waypoints are added to the set of waypoints, denoted as $\mathbb{W}$. At this stage the location of the depots are fixed, and denoted as $\mathbb{D}$. Due to this, the distance function can be represented as a discrete distance matrix, $\mathbb{M}$. The mapping of depot to waypoint can be altered, in the \textit{Assignment Phase}, which is a clear statement of the HTP. The HTP can be solved using the Simplex Algorithm, and alternatively using the Hungarian Algorithm, also known as Kuhn's Combinatorial Algorithm \citep{Ford:1956}. To note further, \cite{Aardal:1998} also investigates the capacitated facility location problem using a polyhedral approach, solved specifically using a cutting-planes approach. However, the implementation of this procedure is beyond the scope of this research due to complexity.

For the \textit{Assignment Phase}, we define a standard assignment problem to be solved using integer programming. We select the Integer Programming approach over the Hungarian Algorithm formulation because the integer programming is a more generalized formulation and does not require redundancy in $\mathbb{M}$.\footnote{See Appendix \ref{appendix} for details.} In the proposed problem, we seek to find a placement of depots, such that the total distance metric from all waypoints to its assigned depot is minimized. We first must define a distance function relating each element the set of waypoints, $w \in \mathbb{W} $, to each element in the set of depots , $d \in \mathbb{D} $. This mapping we define as $\psi$ in Eq. (\ref{eq:f_dist}), 

\begin{equation} \label{eq:f_dist}
    \psi: (\mathbb{W}, \mathbb{D}) \rightarrow \Gamma \in \mathbb{R}
\end{equation}

Where the $\Gamma$ is the set of all distances from from any $w$ to any $d$, which we can also write as $\psi(w, d)$. Thus the distance function $\psi$ is a surjection from the world of n-dimensional coordinate tuples to $\mathbb{R}$. We define this as simply the Euclidean distance between 2-dimensional points in our example. Though it can be extended to the n-dimensional case with the same reasoning. In our use case, since the goal is to minimize the sum of $\Gamma$ \footnote{As defined by the arithmetic sum of all rows and columns of $\Gamma$}. We define $\Delta$, as the assignment plan from each $w$ to $d$ obeying the constraints specified in Eq. (\ref{eq:nk}). $\Delta$ can be represented as a set of tuples $(w, d)$ representing a collection of selected pairings of waypoints to depots. Given this definition our goal is to find a mapping, or assignment, from the set of $\mathbb{W}$ to $\mathbb{D}$ that will be most optimal in terms of minimizing $\psi(w, d)\ \forall \ \Gamma $. We illustrate this relationship where $\mathbb{W} \times \mathbb{D}$ contains all the possible combinations of $w$ and $d$, defined in Eq. (\ref{eq:wxd}).

\begin{align}
    \Delta &\subset \mathbb{W} \times \mathbb{D} \label{eq:wxd} \\ 
    \textrm{s.t.} \quad  \mathbb{D} & \in \mathbb{R}^2 \label{eq:dr} 
\end{align}

Furthermore, we define $\mathbbm{1}_{\Delta}$ as a matrix of indicator variables, stored in a matrix indicating if waypoint $w$ was assigned to depot $d$.

\begin{align} 
\Delta &= 
\begin{bmatrix}
\mathbbm{1}_{\Delta}(w, d)  & ... & \mathbbm{1}_{\Delta}(w, d)\\
\vdots & \ddots & \vdots \\
\mathbbm{1}_{\Delta}(w, d)  & ... & \mathbbm{1}_{\Delta}(w, d)
\end{bmatrix} \label{eq:delta_def}
\end{align}

Our motivation is to find the optimal assignment, denoted as $\Delta^*$, that will minimize the Euclidean distance function . This is illustrated in Eq. (\ref{eq:opt_def}), where $\mathbbm{1}_{\Delta}(w, d) \in (0, 1)$ is an indicator function denoting whether or not $w$ was assigned to $d$ under assignment plan $\Delta$. 

\begin{align}
    \Delta^* &= \underset{\Delta}{\mathrm{argmin}} \sum^{N} \mathbbm{1}_{\Delta}(w, d) \psi(w, d)  \label{eq:opt_def} \\
    &= \underset{\Delta}{\mathrm{argmin}} \ \Psi(\Delta) \label{eq:opt_def_2} 
\end{align}

\section{Solution} \label{sec:solution}

Given a two-phase problem we provide a two-phase algorithm as a solution - outlined in Section (\ref{sec:two_phase}). First we utilize balanced clustering, and subsequently, an integer programming solution to the assignment problem is applied. In Phase I, the balanced clustering algorithm adapts a modified K-Means algorithm to determine the placement of depots freely in $\mathbb{R}^2$, this is known as \textit{balanced clustering}. Subsequently in Phase II, we apply an integer programming solution to the assignment problem, when the depot location have been already determined in Phase I. We provide further details, and an outline, in the Section (\ref{sec:bal_clus}) and (\ref{sec:ass_prob}).

\subsection{Balanced Clustering} \label{sec:bal_clus}

First described in \citep{Macqueen:1967}, the K-Means algorithm presents a geometric interpretation of the classification problem. The algorithm assigns a set of observations into K unconstrained centroids via an iterative algorithm. Consequently we apply this type of algorithm to minimize the total Euclidean distance from each waypoint to its assigned depot. However, in our scenario, we must consider the constraint that $n_k$ must be equal for all groups. For this solution we implement a variation of the k-means algorithms, referred to as the \textit{same-size K-means algorithm} developed by \citep{Schubert:2015}, presented in Algorithm (\ref{algo:elki_equal_km}).

To satisfy the constraint that $n_k$ must be equal for all groups, we constrain the number of waypoints assigned to any $d \in \mathbb{D}$ to be equal. Let $K = |\mathbb{D}|$ denote the cardinality of the set of depots, that is the number of depots. We also constrain the cardinality of each of the assignment subset for any $d$, as denoted by $\Delta_d$, must be equal to the cardinality of all other assignment subsets, represented as $n_k$. We illustrate the simple K-Means algorithm which involves first  assigning each waypoint to its respective depot, which is exactly the closest depot to each waypoint, as illustrated in Eq. (\ref{eq:kmeans_ass}). Subsequently, we re-estimate new depot location, by taking the arithmetic mean of all waypoints assigned to $d$ under $\Delta$, as denoted as $\Delta_d$, as illustrated in Eq. (\ref{eq:kmeans_est}). In K-Means, this alternation between Eq. (\ref{eq:kmeans_ass}) and Eq. (\ref{eq:kmeans_est}) begins initially with a random initialization of candidate depot locations, and ends when either the maximum number of iterations is reached, or when the reduction of $\Psi(\Delta)$ from iteration to iteration is static or below a certain threshold.

\begin{align}
    \Delta(w, d) &= \Big\{ w \ : \ |w - d^*| \leq |w-d|, \ \forall d   \Big\} \label{eq:kmeans_ass} \\
    d^* &= \frac{1}{|\Delta_d|} \sum^{\Delta_d} w \mathbbm{1}_{\Delta}(w, d) \label{eq:kmeans_est} 
\end{align}

Where $\Delta_d$ represents the set of waypoints assigned to depot $d$, and $d^*$ is the proposed optimal depot location. Nevertheless, it is evident that Eq. (\ref{eq:kmeans_ass}) and Eq. (\ref{eq:kmeans_est}) alone does not satisfy the constraint that the cardinality of each $d \in \mathbbm{D}$ to be equal, as illustrated by Eq. (\ref{eq:nk_wd}). Thus to accomplish satisfying Eq. (\ref{eq:nk_wd}), we must apply the algorithm from \citep{Schubert:2015}, and presented in Algorithm (\ref{algo:elki_equal_km}).

\begin{equation} \label{eq:nk_wd}
    |\Delta_d| = \frac{|\mathbbm{W}|}{|\mathbbm{D}|}, \ \forall d \in \mathbbm{D} 
\end{equation}

Using the polynomial time, \textit{Same-size K-Means Algorithm} outlined in Algorithm (\ref{algo:elki_equal_km}), we are able to create a strategy that generates both a set of optimal depot placement locations while maintaining the balanced cluster size constraint from Eq. (\ref{eq:nk_wd}), that is $n_k$ is constant.

\begin{algorithm}[h!]
\caption{Same-size K-Means Algorithm \citep{Schubert:2019}}\label{algo:elki_equal_km}
\begin{algorithmic}[1]
    \State Using K-Means, via Eq. (\ref{eq:kmeans_ass}) and Eq. (\ref{eq:kmeans_est}) propose a set of candidate depot locations.
    \State Compute $\psi(w,d), \ \forall \ \mathbbm{W} \times  \mathbbm{D}$.
    \State Sort $\mathbbm{W}$ based on the difference under $\Delta$ and the best possible alternate assignment, denoted by $\Delta'$, producing ordered set, denoted by $\overline{\mathbbm{W}}$.
    \For {$w \in \overline{\mathbbm{W}} $}:
        \For {$d \in \ \mathbbm{D} $}:
            \State Initialize $\overline{\mathbbm{W}}_{d'}$ as all $w \in \Delta_{d'}$
                \While {$|\overline{\mathbbm{W}}_{d'}| > 0$}:
                    \For {$w$ in $\Delta_d$}
                        \If {Swapping $w$ with $w'$ from $d$ to $d'$ reduces $ \Psi(\Delta^{'}) $}:
                            \State Assign $w$ to $d'$ and $w'$ to $d$.
                            \State Remove $w$ from $\overline{\mathbbm{W}}_d$.
                        \EndIf
                        \If {Reassigning $w$ to $d'$ does not violate Eq. (\ref{eq:nk_wd})}:
                            \State Assign $w$ to $d'$.
                            \State Remove $w$ from $\overline{\mathbbm{W}}_d$.
                        \EndIf
                        \If {Maximum iterations reached}
                            \State Terminate algorithm.
                        \EndIf
                    \EndFor
                \EndWhile
        \EndFor
    \EndFor
\end{algorithmic}
\end{algorithm}

\subsection{The Assignment Problem} \label{sec:ass_prob}

After the depot locations have been determined in Phase I, as outlined previously in Section \ref{sec:bal_clus}, we can formulate a tractable solution for the \textit{Assignment Phase}. This can be constructed as an Hitchcock Transportation Problem. Eq. (\ref{eq:ass_def_depot}) illustrates the objective function and constraints for such an optimization problem. In our formulation, we assign to each waypoint, $w$ to a fixed depot, $d$, where $\mathbbm{D}$ no longer lies unconstrained in $\mathbbm{R}^2$. $\mathbbm{D}$ is a fixed set, which we denote as $\widetilde{\mathbbm{D}}$. Provided $\widetilde{\mathbbm{D}}$, we introduce a set of new fixed waypoints $\widetilde{\mathbbm{W}}$, where we must assign each $w \in \widetilde{\mathbbm{W}}$ to a specific depot in $\widetilde{\mathbbm{D}}$. Alternatively, \cite{Malinen:2014} suggests that using the Hungarian Algorithm, also called Munkres algorithm, \citep{Kuhn:1955} is capable to solve up to approximately 1000 waypoints on regular computers for this specific problem. However, in order to formulate the problem using the Hungarian Algorithm, it is necessary to repeat the depot locations on the cost matrix - see Appendix \ref{appendix}. This is inefficient, and we opt for a direct Integer Programming (IP) solution which we will present in this paper. We present a solution where we optimally assign each waypoint to depot such that $\Psi(\Delta)$ is minimized.

\begin{align} 
    \min_{ \Delta } \quad & \sum_{ \Gamma }{\mathbbm{1}_{\Delta}(w, d)\psi(w, d)  } \label{eq:ass_goal_func} \\
    \textrm{s.t.} \quad & \sum_{\widetilde{\mathbbm{D}}}{ \mathbbm{1}_{\Delta}(w, d) } = \frac{|\widehat{\mathbbm{W}}|}{|\widetilde{\mathbbm{D}}|}, \ \forall w \in \widehat{\mathbbm{W}}  \label{eq:ass_def_wp} \\
    & \sum_{\widehat{\mathbbm{W}}}{ \mathbbm{1}_{\Delta}(w, d) } = 1, \ \forall d \in \widetilde{\mathbbm{D}}  \label{eq:ass_def_depot}
\end{align}

In Eq. (\ref{eq:ass_goal_func}) we illustrate the objective function that must be minimized in our IP, and subsequently Eq. (\ref{eq:ass_def_wp}) and (\ref{eq:ass_def_depot}) we illustrate the constraints on the IP. In our experiment, previously assigned waypoints assigned in Phase I may be reassigned to another depot, however, the depot locations $\widetilde{\mathbbm{D}}$ are fixed. We use $\widehat{\mathbbm{W}}$ to denote the final set of waypoints, as marked by Eq. (\ref{eq:w_union}), as the union of the waypoints issued in both Phase I and Phase II. The assignment plan $\Delta$ can be constructed as a matrix of indicator variables $\mathbbm{1}_{\Delta}(w, d)$, indicating whether waypoint $w$ was assigned to depot $d$, as illustrated in Eq. (\ref{eq:delta}). As evident from Eq. (\ref{eq:wayp_const}) and (\ref{eq:depot_const}), we specify the constraints of the way points as the assignment plan $\Delta$ and the transpose of the assignment plan $\Delta^T$. Each row of the assignment plan $\Delta$ refers to a specific waypoint, $w$, and each column, a depot $d$. As a small non-convention, but for the sake of clear simplicity we define the row-sum operator $[\Delta]^+$ a new vector containing the matrix row sums of the matrix $\Delta$. Allowing us to clearly specify the constraints.

\begin{align} 
\Delta &= 
\begin{bmatrix}
\mathbbm{1}_{\Delta}(w, d) \ & ... & \ \mathbbm{1}_{\Delta}(w, d)\\
\vdots & \ddots & \vdots \\
\mathbbm{1}_{\Delta}(w, d) \ & ... & \ \mathbbm{1}_{\Delta}(w, d)
\end{bmatrix} \label{eq:delta} \\ \
[\Delta]^+ &= 
\begin{bmatrix}
\mathbbm{1}_{\Delta}(w, d) \ + & ... & + \ \mathbbm{1}_{\Delta}(w, d)\\
\vdots & \ddots & \vdots \\
\mathbbm{1}_{\Delta}(w, d) \ + & ... & + \ \mathbbm{1}_{\Delta}(w, d)
\end{bmatrix} =
\begin{bmatrix}
1 \\
\vdots \\
1 \label{eq:wayp_const} 
\end{bmatrix} \\ \
[\Delta^T]^+ &=
\begin{bmatrix}
n_k & ... & n_k
\end{bmatrix}^T \label{eq:depot_const} 
\end{align}

Several industry standard solutions to the IP problem exist for the NP-hard Hitchcock Problem. We present some preliminary background on such methods, however, we will not go into heavy theoretical details, as it is not the focus of our thesis. The solutions to IP problems concern the Cutting Planes method \citep{Gomory:1958}, which includes a solution to the separation problem, that is, finding an inequality that separates the optimal value from the convex hull, known as a \textit{cut}. Many cuts are found until the non-integer optimal solution is no longer feasible. The cutting plane method, was in its introduced in the 1950's had impractical applications due to numerical instabilities. Another solution to the IP problem is known as Branch-and-Bound \citep{LandDoig:1960}. In principle, the Branch-and-Bound algorithm involves iterative and selective computations to produce subspaces of the feasible solution after elimination and eliminating such subspaces by remembering the bounds of the solution space pertaining so such a subspace, and eliminating them accordingly. This algorithm expedites the computation of the iterative search involved in solving for the global optimal.

We option to apply the Branch-and-Cut \citep{Balas:1991} method as the solution to our IP, on the grounds that it combines the effectiveness of both the Cutting Planes and Branch-and-Bound methods, improving the numerical instability of the Cutting planes method, while also improving the subspace elimination capabilities of the original Branch-and-Bound. This serves currently, as a state-of-the-art solution for most Mixed Integer Programs today. Therefore, provided the objective functions and constraints illustrated in Eq. (\ref{eq:ass_goal_func})
, (\ref{eq:wayp_const}), and (\ref{eq:wayp_const})., we proceed to solve the system of equations using the Branch-and-Cut method \citep{Elf:2001}. We use the software implementation from the open source package for Branch-and-Cut developed by COIN-OR, \textit{Computational Infrastructure for Operations Research} \citep{coinor:2020}. 

\subsection{Two-Phase Algorithm}\label{sec:two_phase}

The full algorithm proposed in this work begins with a set of initial waypoints, $\mathbbm{W}$, and then we are subsequently introduced new waypoints, $\widetilde{\mathbbm{W}}$. As outlined, we define two phases of optimization, the balanced clustering, and Assignment Phases. In both phases, the depot locations, $\mathbbm{D}$ remains constant, but are positioned in Phase I. In reality, $\mathbbm{D}$ can  represent the required placement of distribution centers that must be placed to serve such locations. The initial set of waypoints, $\mathbbm{W}$, can represent, for example, store locations or service centers, that are initially planned by a company or organization. The second phase of waypoints, $\widetilde{\mathbbm{W}}$, can represent a subsequent addition of store locations planned by a company. Typically, the planning of these waypoints are perhaps strategic, however, we assume the locations of these stores to be unknown in Phase I. In our expriment, we simulate this strategic placement of new stores, by randomly sampling from a set of $\widehat{\mathbbm{W}}$ all store locations, randomly allocating $\mathbbm{W}$ to the balanced clustering phase, and $\widetilde{\mathbbm{W}}$ to the Assignment phase, with the percentage of allocation respecting $\gamma = 0.05$, as described in Section \ref{sec:data_desc}.

\begin{equation} \label{eq:w_union}
    \widehat{\mathbbm{W}} = \mathbbm{W} \cup \widetilde{\mathbbm{W}}
\end{equation}

The balanced clustering phase utilizes the iterative algorithm outlined in Algorithm (\ref{algo:elki_equal_km}) to generate a prescribed set of depot locations. In this phase, the depot locations fall anywhere on an $\mathbbm{R}^2$ plane. Inevitably, the assignment, $\Delta$, is simply the closest depot, $d$, to each $w \in \mathbbm{W}$. In the second phase, we perform the assignment algorithm illustrated in Section (\ref{sec:ass_prob}). Where, as described, we assign new waypoints, $\widetilde{\mathbbm{W}}$ to $\mathbbm{D}$. The second phase can represent new stores or service locations that are planned for the future. Since the depots have been already constructed, we are only allowed to solve the assignment problem with $\psi(w,d)$ serving as the loss measure, of which we minimize. Algorithm \ref{algo:two_phase} summarizes this Two-Phase strategy. We also provide source code for implmenting this algorithm on \citep{Liu_git:2020}. 

\begin{algorithm}[h!]
\caption{Two-Phase Algorithm}\label{algo:two_phase}
\begin{algorithmic}[1]
    \State $\mathbbm{W}$ is given.
    \State Compute $\Delta^*$ using Algorithm (\ref{algo:elki_equal_km})
    \State $\widehat{\mathbbm{W}}$ is given.
    \State Compute $\Delta^*$ using Integer Programming solution outlined in Section (\ref{sec:ass_prob}).
\end{algorithmic}
\end{algorithm}

% \begin{figure}[!h]
%   \includegraphics[width=\textwidth]{Balanced_Clustering_Flow.png}
%   \caption{balanced clustering Process}
%   \label{fig:Bal_Clus_Flow}
% \end{figure}

% \begin{figure}[!h] 
%   \includegraphics[width=\textwidth]{9_depot_sb_map.png}
%   \caption{Optimal placement of depot locations for 9 service depots around the world.}
%   \label{fig:depot_locs}
% \end{figure} 

\begin{figure}[!h] 
  \includegraphics[width=\textwidth]{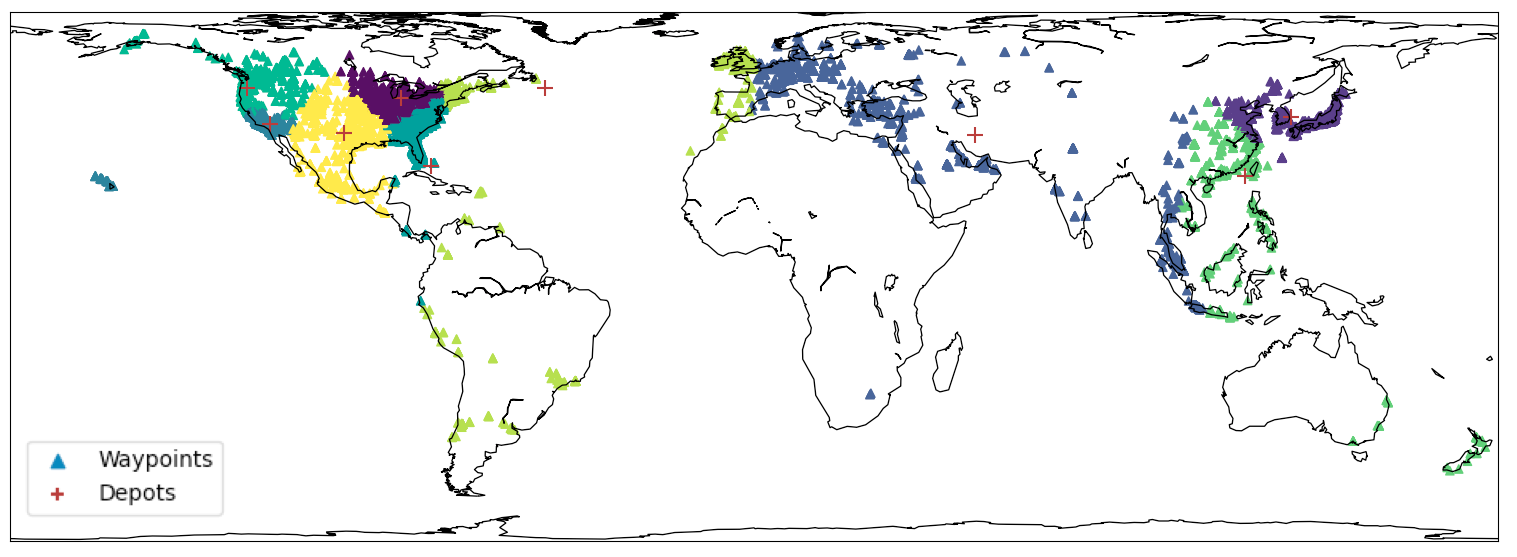}
  \caption{Optimal placement of depot locations for K = 9 service depots around the world.}
  \label{fig:ass_color}
\end{figure} 

\section{Simulation}

In order to measure the performance of the Two-Phase algorithm, we use the Mean Squared Error (MSE) of distance function $\psi$ of each $w \in \widehat{\mathbbm{W}}$ to its assigned depot, as designated as $\Delta^*[w]$.

\begin{equation} \label{eq:mse}
     MSE = \frac{1}{|\widehat{\mathbbm{W}}|} \sum_{w \in \widehat{\mathbbm{W}}} {\psi(w, \Delta^*[w] ) ^2} 
\end{equation}

Eq. (\ref{eq:mse}) produces a measure of the average of cost from the each waypoint, $w$ to its depot mapping $\Delta^*[w]$. This is a common metric for measuring clustering performance, as we wish to reduce the Euclidean distance from $w$ to its optimally mapped depot, $d = \Delta^*[w]$. As evidence from the performance from Table (\ref{table:stat_results}), we see that the percent deviation from Phase I to Phase II in terms of MSE is under on average of 3\%, and no greater than 10 \% for up $K = 10$ starting from $K = 3$. This indicates that, even if mobility of depots is not allowed in Phase II, the Assignment Algorithm using Integer Programming still produces a mapping $\Delta$ that does not differ far from a supposedly more flexible, and thereby more optimal Phase I algorithm. In general, this strategy is simple to understand and implement, and provides strong empirical results.

\begin{table}[h!]
\centering
\begin{tabular}{||c c c c||} 
 \hline
 K & MSE I & MSE II & \% Change \\ [0.5ex] 
 \hline\hline
 3&985.428 & 987.747 & 0.002341\\ 
 4&935.636 & 945.533 &  0.010467\\ 
 5&903.949 & 920.728 & 0.018224\\ 
6&651.008&655.248 & 0.006470 \\ 
7&470.978&480.461 & 0.019737 \\ 
8&298.556&332.964 & 0.103338 \\ 
9&450.548&459.597 & 0.019688 \\ 
10&348.791&373.617 & 0.066447 \\ [1ex] 
 \hline
\end{tabular}
\caption{Comparison between MSE Phase I, and MSE Phase II. }
\label{table:stat_results}
\end{table}

\section{Conclusion}

In our research we present a basic logistics problem, a two-phase facility planning problem, with a different set of constraints in each phase. In each of the two phases, we present an application of pre-existing methods to solve it. Notably, the \textit{same size K-means} algorithm \citep{Schubert:2015} for the balanced clustering phase, and an Integer Programming solution to the Hitchcock Transportation problem utilizing the Branch-and-Cut algorithm. Furthermore, we present a unique mathematical formulation that ties the two phases together in a unified notation. Subsequently, we prove the capability of our proposed two phase algorithm in a simulation framework on real world data. We acknowledge that further research can be done with relation to some of the hyperparameters of the project. Hyperparameters include the amount of waypoints given at each phase, as determined by $\gamma$, and also for example, the maximum number of iterations of the \textit{same size K-means} algorithm. We can also potentially study the efficacy of this two-phase algorithm beyond two dimensions. Or also, increasing the number of depots, K, and of course, attempting to run such an algorithm on other datasets. Nevertheless, this paper serves as a proposal for a possibility of applying such a two-phase algorithm any multi-phase planning of facilities as illustrated in our work.

\section*{Acknowledgements}

This research was conducted without external funding, resulting in no conflicts of interests. The sole author would like to thank the additional efforts of Dr. Seyed Ali Hesammohseni from the University of Waterloo for providing extensive comments regarding the scientific direction of the work, as well as the Kaggle oganization for providing the open source data set.

\appendix
\section{} \label{appendix}

% \subsection{Summary of Notation Definitions}\label{app:notation}

% \subsection{Note on the Hungarian Algorithm} \label{app:hung}

We illustrate a problem statement for alternate solution to the Hitchcock Transportation problem, using the Hungarian Method \citep{Kuhn:1955}. This method was purposed to find a minimum cost solution for assignment of workers to tasks, given the capacity of the workers, and the requirements for the each of the tasks. Notably, the original Hungarian Algorithm requires that the number of workers equals the number of tasks, where as in our scenario the number of waypoints greatly exceed the number of depots. It applies a series of algebraic manipulations to compute an optimal solution for the Assignment Problem. This involves only a one worker to one task solution based on a matrix that contains the cost per worker per task. In order to address this mismatch, we can duplicate the number of depots by $n_k$ times, increasing the number of columns repetitively, to create a $N$ by $K \cdot n_k$ square matrix. Let $\mathbbm{M}$ represent the distance matrix from each $w \in \mathbbm{W}$ to $d \in \mathbbm{D}$. Therefore we have an $N \times K$ distance matrix, which we denote as $\widetilde{\mathbbm{M}}$. 

\begin{equation} \label{eq:M_mat}
     \mathbbm{M} = \begin{bmatrix}
\psi(w, d) & ... & \psi(w, d)\\
\vdots & \ddots & \vdots \\
\psi(w, d) & ... & \psi(w, d)
\end{bmatrix}
\end{equation}

In order to build a matrix where the Hungarian can be applied we must construct $\widetilde{\mathbbm{M}}$, we simply replicate $\mathbbm{M}$ column-wise a total of $n_k$ times. In this respect, we can obtain an $N \times K$ matrix. 

\begin{equation} \label{eq:M_tilde_mat}
     \widetilde{\mathbbm{M}} = [\mathbbm{M} \ ... \ \mathbbm{M} ]
\end{equation}

Because of the repetition created in Eq. (\ref{eq:M_tilde_mat}), we opt to not apply the Hungarian as a viable solution for our experiment when $\Delta^*$ in the Assignment Phase.

\vskip 0.2in
\bibliography{sample}

\begin{thebibliography}{19}
\providecommand{\natexlab}[1]{#1}
\providecommand{\url}[1]{\texttt{#1}}
\expandafter\ifx\csname urlstyle\endcsname\relax
  \providecommand{\doi}[1]{doi: #1}\else
  \providecommand{\doi}{doi: \begingroup \urlstyle{rm}\Url}\fi

\bibitem[Aardal(1998)]{Aardal:1998}
Karen Aardal.
\newblock Capacitated facility location: Separation algorithms and
  computational experience.
\newblock \emph{Math. Program.}, 81:\penalty0 149--175, 1998.
\newblock \doi{10.1007/BF01581103}.
\newblock URL \url{https://doi.org/10.1007/BF01581103}.

\bibitem[Balas and Matthew J.~Saltzman(1991)]{Balas:1991}
Egon Balas and Gérard~Cornuéjols Matthew J.~Saltzman.
\newblock An algorithm for the three-index assignment problem.
\newblock \emph{Oper. Res.}, 39\penalty0 (1):\penalty0 150--161, 1991.
\newblock \doi{10.1287/opre.39.1.150}.
\newblock URL \url{https://doi.org/10.1287/opre.39.1.150}.

\bibitem[Dempster et~al.(1977)Dempster, Laird, and Rubin]{Dempster:1977}
A.~P. Dempster, N.~M. Laird, and D.~B. Rubin.
\newblock Maximum likelihood from incomplete data via the {EM} algorithm.
\newblock \emph{Journal of the Royal Statistical Society: Series B},
  39:\penalty0 1--38, 1977.
\newblock URL \url{http://web.mit.edu/6.435/www/Dempster77.pdf}.

\bibitem[Elf et~al.(2001)Elf, Gutwenger, J{\"u}nger, and Rinaldi]{Elf:2001}
Matthias Elf, Carsten Gutwenger, Michael J{\"u}nger, and Giovanni Rinaldi.
\newblock \emph{Branch-and-Cut Algorithms for Combinatorial Optimization and
  Their Implementation in ABACUS}, pages 157--222.
\newblock Springer Berlin Heidelberg, Berlin, Heidelberg, 2001.
\newblock ISBN 978-3-540-45586-8.
\newblock \doi{10.1007/3-540-45586-8_5}.
\newblock URL \url{https://doi.org/10.1007/3-540-45586-8_5}.

\bibitem[Fekete et~al.(2003)Fekete, Mitchell, and Beurer]{Fekete:2003}
S{\'{a}}ndor~P. Fekete, Joseph S.~B. Mitchell, and Karin Beurer.
\newblock On the continuous fermat-weber problem.
\newblock \emph{CoRR}, cs.CG/0310027, 2003.
\newblock URL \url{http://arxiv.org/abs/cs/0310027}.

\bibitem[Ford and Fulkerson(1956)]{Ford:1956}
L.~R. Ford and D.~R. Fulkerson.
\newblock Solving the transportation problem.
\newblock \emph{Management Science}, 3\penalty0 (1):\penalty0 24--32, 1956.
\newblock ISSN 00251909, 15265501.
\newblock URL \url{http://www.jstor.org/stable/2627172}.

\bibitem[Forest et~al.(2020)Forest, Vigerske, Ralphs, Hafer, jpfasano, Santos,
  Saltzman, h-i gassmann, Kristjansson, and King]{coinor:2020}
Jon Forest, Stefan Vigerske, Ted Ralphs, Lou Hafer, jpfasano, Haroldo~Gambini
  Santos, Matthew Saltzman, h-i gassmann, Bjarni Kristjansson, and Alan King.
\newblock coin-or/clp: Version 1.17.6, April 2020.
\newblock URL \url{https://doi.org/10.5281/zenodo.3748677}.

\bibitem[Gomory(1958)]{Gomory:1958}
Ralph~E. Gomory.
\newblock Outline of an algorithm for integer solutions to linear program.
\newblock \emph{Bulletin of the American Mathematical Society}, 64\penalty0
  (5):\penalty0 275--278, September 1958.

\bibitem[Hakimi(1964)]{Hakimi:1964}
S.~L. Hakimi.
\newblock Optimum locations of switching centers and the absolute centers and
  medians of a graph.
\newblock \emph{Operations Research}, 12\penalty0 (3):\penalty0 450--459, 1964.
\newblock URL
  \url{https://EconPapers.repec.org/RePEc:inm:oropre:v:12:y:1964:i:3:p:450-459}.

\bibitem[Hitchcock(1941)]{Hitchcock:1941}
Frank~L. Hitchcock.
\newblock The distribution of a product from several sources to numerous
  localities.
\newblock \emph{Journal of Mathematics and Physics}, 20\penalty0
  (1-4):\penalty0 224--230, 1941.
\newblock \doi{10.1002/sapm1941201224}.
\newblock URL
  \url{https://onlinelibrary.wiley.com/doi/abs/10.1002/sapm1941201224}.

\bibitem[Kuhn(1955)]{Kuhn:1955}
Harold~W. Kuhn.
\newblock The hungarian method for the assignment problem.
\newblock \emph{Naval Research Logistics Quarterly}, 2:\penalty0 83--97, 1955.

\bibitem[Land and Doig(1960)]{LandDoig:1960}
A.~H. Land and A.~G. Doig.
\newblock An automatic method of solving discrete programming problems.
\newblock \emph{Econometrica}, 28\penalty0 (3):\penalty0 pp. 497--520, 1960.
\newblock ISSN 00129682.

\bibitem[Liu(2020)]{Liu_git:2020}
Larkin Liu.
\newblock Two-phase algorithm for facility planning.
\newblock \url{https://github.com/larkz/equal_clustering}, 2020.

\bibitem[MacQueen(1967)]{Macqueen:1967}
J.~B. MacQueen.
\newblock Some methods for classification and analysis of multivariate
  observations.
\newblock In L.~M.~Le Cam and J.~Neyman, editors, \emph{Proc. of the fifth
  Berkeley Symposium on Mathematical Statistics and Probability}, volume~1,
  pages 281--297. University of California Press, 1967.

\bibitem[Malinen and Fr{\"a}nti(2014)]{Malinen:2014}
Mikko~I. Malinen and Pasi Fr{\"a}nti.
\newblock Balanced k-means for clustering.
\newblock In Pasi Fr{\"a}nti, Gavin Brown, Marco Loog, Francisco Escolano, and
  Marcello Pelillo, editors, \emph{Structural, Syntactic, and Statistical
  Pattern Recognition}, pages 32--41, Berlin, Heidelberg, 2014. Springer Berlin
  Heidelberg.
\newblock ISBN 978-3-662-44415-3.

\bibitem[Schubert and Zimek(2019)]{Schubert:2019}
Erich Schubert and Arthur Zimek.
\newblock {ELKI:} {A} large open-source library for data analysis - {ELKI}
  release 0.7.5 "heidelberg".
\newblock \emph{CoRR}, abs/1902.03616, 2019.
\newblock URL \url{http://arxiv.org/abs/1902.03616}.

\bibitem[Schubert et~al.(2015)Schubert, Koos, Emrich, Z{\"{u}}fle, Schmid, and
  Zimek]{Schubert:2015}
Erich Schubert, Alexander Koos, Tobias Emrich, Andreas Z{\"{u}}fle,
  Klaus~Arthur Schmid, and Arthur Zimek.
\newblock A framework for clustering uncertain data.
\newblock \emph{Proc. {VLDB} Endow.}, 8\penalty0 (12):\penalty0 1976--1979,
  2015.
\newblock \doi{10.14778/2824032.2824115}.
\newblock URL \url{http://www.vldb.org/pvldb/vol8/p1976-schubert.pdf}.

\bibitem[Voronoi(1908)]{Voronoi:1908}
Grigory Voronoi.
\newblock Nouvelles applications des paramètres continus à la théorie des
  formes quadratiques. premier mémoire. sur quelques propriétés des formes
  quadratiques positives parfaites.
\newblock \emph{Journal für die reine und angewandte Mathematik}, 1908.
\newblock URL
  \url{https://www.deutsche-digitale-bibliothek.de/item/KDI2RGXOO22HGVX6KZCPNRI52ZTZCHVG}.

\bibitem[Weber(1909)]{Weber:1909}
A.~Weber.
\newblock \emph{Ueber den Standort der Industrien: Reine Theorie des Standorts,
  mit einem mathematischen Anhang, von Georg Pick}.
\newblock Ueber den Standort der Industrien. J.C.B. Mohr (Paul Siebeck), 1909.
\newblock URL \url{https://books.google.de/books?id=FSrZAAAAMAAJ}.

\end{thebibliography}

\end{document}